\documentclass[11pt]{article}

\usepackage{natbib}
 \bibpunct[, ]{(}{)}{,}{a}{}{,}%

\usepackage{multirow}
\usepackage{amsmath,amssymb}
\usepackage[ruled, vlined]{algorithm2e}
\usepackage{graphicx}
\usepackage{tikz}
\usepackage{booktabs} 
\usepackage{xcolor}
\usepackage{colortbl}
\usepackage{wrapfig}

\DeclareGraphicsExtensions{.pdf,.png,.jpg}

\def\nexto{\kern -0.54em}

\def\DoubleZ{{\cal Z \kern -0.7em Z}}

\newcommand{\bp}{\mathcal{P}}
\newcommand{\bq}{\mathcal{Q}}
\newcommand{\bpq}{\mathcal{QP}}

\newtheorem{observation}{Observation}

\begin{document}

\begin{center}

\LARGE

{\bf Stratified Patient Appointment Scheduling \\
for Community-based \\
Chronic Disease Management Programs \\}

\normalsize

\vspace{24pt}

{Martin Savelsbergh} \\
{\it Georgia Institue of Technology, Atlanta} \\
\vspace{12pt}
{Karen Smilowitz} \\
{\it Northwestern University, Evanston}
\end{center}

\begin{abstract}

Disease management programs have emerged as a cost-effective approach to treat chronic diseases. Appointment adherence is critical to the success of such programs; missed appointment are costly, resulting in reduced resource utilization and worsening of patients' health states. The time of an appointment is one of the factors that impacts adherence. We investigate the benefits, in terms of improved adherence, of incorporating patients' time-of-day preferences during appointment schedule creation and, thus, ultimately, on population health outcomes. Through an extensive computational study, we demonstrate, more generally, the usefulness of patient stratification in appointment scheduling in the environment that motivates our research, an asthma management program offered in Chicago. We find that capturing patient characteristics in appointment scheduling, especially their time preferences, leads to substantial improvements in community health outcomes. We also identify settings in which simple, easy-to-use policies can produce schedules that are comparable in quality to those obtained with an optimization-based approach.
\end{abstract}

\noindent \textbf{Keywords:} appointment scheduling, chronic
disease, community-based care, disease progression, patient no-show, time-of-day preference

\section{Introduction}

Disease management programs have emerged as a cost-effective
approach to treat chronic diseases; see \cite{jones2007achieving}. A
disease management program serves a patient population for a
specific chronic disease, such as asthma or diabetes; see
\cite{jones2005breathmobile} and \cite{kucukyazici2013managing}. The
asthma management program offered in Chicago by The Mobile C.A.R.E.
Foundation (MCF) is an example of such a disease management program.
Through a partnership with the Chicago public school system, MCF
serves asthmatic children through repeated visits to their school
with mobile clinics.  At each visit, asthmatic children are examined
and treated by a medical team.  The visits, as well as the children
that will be seen during the visits, are scheduled months in advance
by an MCF administrator.

The characteristics of this setting lead to interesting challenges in appointment scheduling.  Firstly, community-based disease management programs (such as the program that motivates our research) typically serve a fixed patient population with limited care capacity and with a goal of maximizing health outcomes for the entire population.  Secondly, the nature of chronic conditions requires recurring patient visits over a planning horizon to
maintain disease control. Disease progression occurs between visits,
which must be taken into account when scheduling appointments.
Thirdly, unlike traditional appointment scheduling settings in which
patients request appointments as needed, appointments in
community-based chronic disease management programs are scheduled by
the provider.  This often happens far in advance and overbooking
appointment slots may not be possible. In the case of MCF, for
example, privacy issues and a lack of space to wait in the mobile
clinic prohibit overbooking. Finally, because appointments are
scheduled far in advance, there is a higher likelihood that patients
fail to show up for an appointment.
\begin{wrapfigure}{r}{0.4\textwidth}
	\includegraphics[scale=0.4]{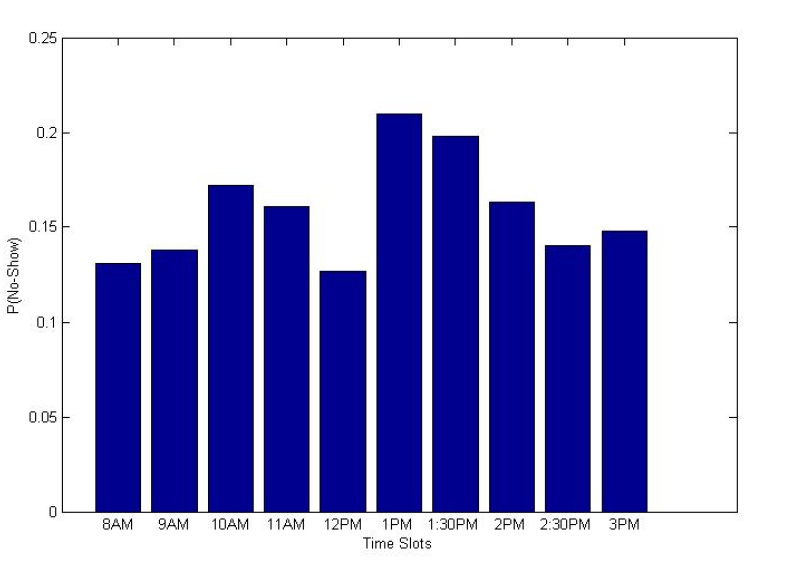}
	\vspace{-1cm}
	\caption{\small Historical missed appointment percentages by time of day at
		MCF. \label{f:mcf-missed}}
\end{wrapfigure}
In the case of MCF, a parent or guardian must accompany the patient, and no-show rates of more than 15\% are not uncommon. Importantly, an analysis of historical data from MCF shows that no-show rates vary with time of day, with lower no-show rates in the early morning, during lunch, and the late afternoon; see Figure \ref{f:mcf-missed}. This is due mostly to the work schedules of parents who must accompany a patient. Consequently, it appears to be important to consider not only the interval between visits, but also the time of the visits when
creating appointment schedules. This explains the two main goals of our research: (1) to assess the benefits, in terms of health outcomes, of accounting for time-of-day preferences in appointment scheduling, and (2) to investigate how easy or hard it is to incorporate time-of-day preferences in appointment scheduling procedures.

We explore both a sophisticated appointment scheduling method, which
considers individual patient characteristics, and relatively simple
and easy-to-use appointment scheduling methods, which only
distinguish groups of patients with similar characteristics, which
we refer to as ``cohort scheduling policies''.  These approaches are
compared based on their ability to maximize the health
state of a population, measured by the likelihood that patients'
disease is controlled.  We design a set of stylized test instances
based on our motivating setting to better understand the importance
of accounting for patient-specific, time-dependent no-show rates.
The study demonstrates that explicitly accounting for these factors
produces appointment schedules with substantially better population
health outcomes, up to 15\% better in some settings. The study
further shows that easy-to-use cohort-based methods are effective in
settings with a fairly homogeneous patient population and in settings
in which patient preferences are known or can easily be deduced.
These results are encouraging and highlight the tremendous potential
of acquiring and using patients' time-of-day preferences to
construct more effective appointment schedules resulting in better
population health outcomes.

The remainder of the paper is organized as follows. In Section
\ref{sec:char}, we present relevant literature and discuss the
characteristics of the appointment scheduling environment we
consider. In Section \ref{s:sch_app}, we present optimization-based
and cohort-based appointment scheduling approaches. In Section
\ref{s:comp_study}, we present a computational study, in which we
evaluate the performance of  scheduling methods for patient
populations with different characteristics. We end with final
remarks in Section \ref{s:final}.

\section{Scheduling patient appointments in chronic disease management}
\label{sec:char}

\cite{kucukyazici2013managing} present an overview of
community-based care programs for chronic diseases, detailing
program operations and effectiveness and highlighting three relevant
papers from the operations research literature: \cite{leff1986lp},
\cite{Deoetal2011}, and \cite{kucukyazici2011analytical}.  Critical
to all of these papers is the interaction between care provided and
patient health state, as patients' health states change over time
and with access to health care.  Our work builds on
\cite{Deoetal2011}, which first examined the challenges of
appointment scheduling for MCF. In that paper, the authors present
an integrated capacity allocation model to select which patients to
see each period.  The model combines clinical (disease progression) and
operational (capacity constraint) factors and is shown to outperform
traditional strategies that decouple the two. However, the model
does not consider patient no-shows, and, as a result, the allocation
of time slots within a day is not part of the model.

Because overbooking is not an option for MCF, each patient no-show
implies a loss of already scarce provider capacity. \cite{MCFcase09}
highlights the impact of patient no-show rates on MCF operations and
outlines the strict policies in place to ensure that provider
capacity is used as effectively as possible: patients who miss
appointments repeatedly may be removed from the program, and schools
with excessive aggregate no-show rates may also be removed from the
program. In this paper, by explicitly considering the temporal
dependence of patient no-show rates in appointment scheduling, we
hope to reduce the occurrence of no-shows.

Appointment scheduling problems in health care settings have been
the focus of much recent work.  \cite{gupta2008appointment} provide
a comprehensive overview of the subject; more recent surveys have
included appointment scheduling in reviews of operations research
techniques in a wider range of health care decisions; see
\cite{batun2013optimization} and \cite{hulshof2012taxonomic}.
\cite{gupta2008appointment} categorize the health care appointment
scheduling literature by scheduling environment, given the unique
characteristics of each setting: primary care, speciality clinic,
and elective surgery.  Our setting, i.e., a chronic disease
management program, shares some characteristics with speciality
clinics and elective surgery, with a few key differences.  As with
elective surgery, patients are scheduled in a ``single batch'',
meaning an administrator schedules all slots for a given time period
at once, rather than scheduling appointments as patients make
requests, as is the case in primary care and often in speciality
clinics.  However, unlike surgical settings, chronic patients
require recurring visits to the provider and the interval between
these visits impacts disease progression; see
\cite{jones2007achieving}.

\cite{SSV2008} demonstrate the relevance of appointment adherence in
a study of the impact of no-shows among patients with diabetes.  The
authors find that for each 10\% increment in missed appointment
rate, the odds of good control decrease by a factor of 1.12 and the
odds of poor control increase by factor of 1.24.
\cite{gupta2008appointment} identify patient no-shows as a key
factor in scheduling and highlight approaches, e.g., open access and
overbooking, to address no-shows; see \cite{robinson2010comparison}
and \cite{liu2010dynamic} for examples of recent work. However, as
noted earlier, the lack of waiting space and the requirement that
parents accompany patients preclude such options in our setting.

A few recent papers have considered patient preferences (including
time-of-day preferences) in the presence of no-shows; see
\cite{feldman2012appointment}, \cite{gupta2008revenue}, and
\cite{wang2011adaptive}. These papers  consider dynamic
settings in which patients are scheduled as they make requests.
\cite{feldman2012appointment} consider a multi-day setting in which
patients are offered a set of appointment options at the time of the
appointment request.  While their paper also considers a static
model, the request arrival and scheduling setting are quite different
from the single batch scheduling in our chronic disease management
program setting.  Their work, along with the others in this stream,
are more akin to revenue management models.

\cite{samoranioutpatient} study multi-day dynamic appointment
scheduling with no-show rates that vary by patient and time since
booking in an outpatient setting with open access and overbooking.
Their work includes a detailed analysis of data from a mental health
center to identify causes of no-shows among patients.  Unlike our
setting, appointments are scheduled on a rolling horizon as patients
request appointments.  While these differences lead to a
fundamentally different model, the authors use
column generation to handle the large number of variables in their
model, similar to our optimization-based approach.

Next, we present the key features and characteristics of the patient
appointment scheduling environment we consider: patient appointment
schedule, patient disease progression, patient no-show
probabilities, and patient appointment schedule evaluation.

\subsection{Patient appointment schedule} \label{sb:pat_app_sch}

We consider an appointment scheduling environment where the planning horizon consists of $K$ periods, each with $T$ time slots, and where there are $P$ patients in the population. In the context of MCF, a period represents a day in which a mobile clinic visits a particular school with $P$ asthmatic students, $T$ of which can be seen that day. For ease of notation and consistent with MCF practice, we assume that periods are equally spaced in time; however, our models can be generalized if this is not the case. Due to the limited number of time slots in each period, it is typically not possible to see all patients each period (i.e., $P > T$).

\begin{figure}[htbp]
    \centering
        \includegraphics[width=0.60\textwidth]{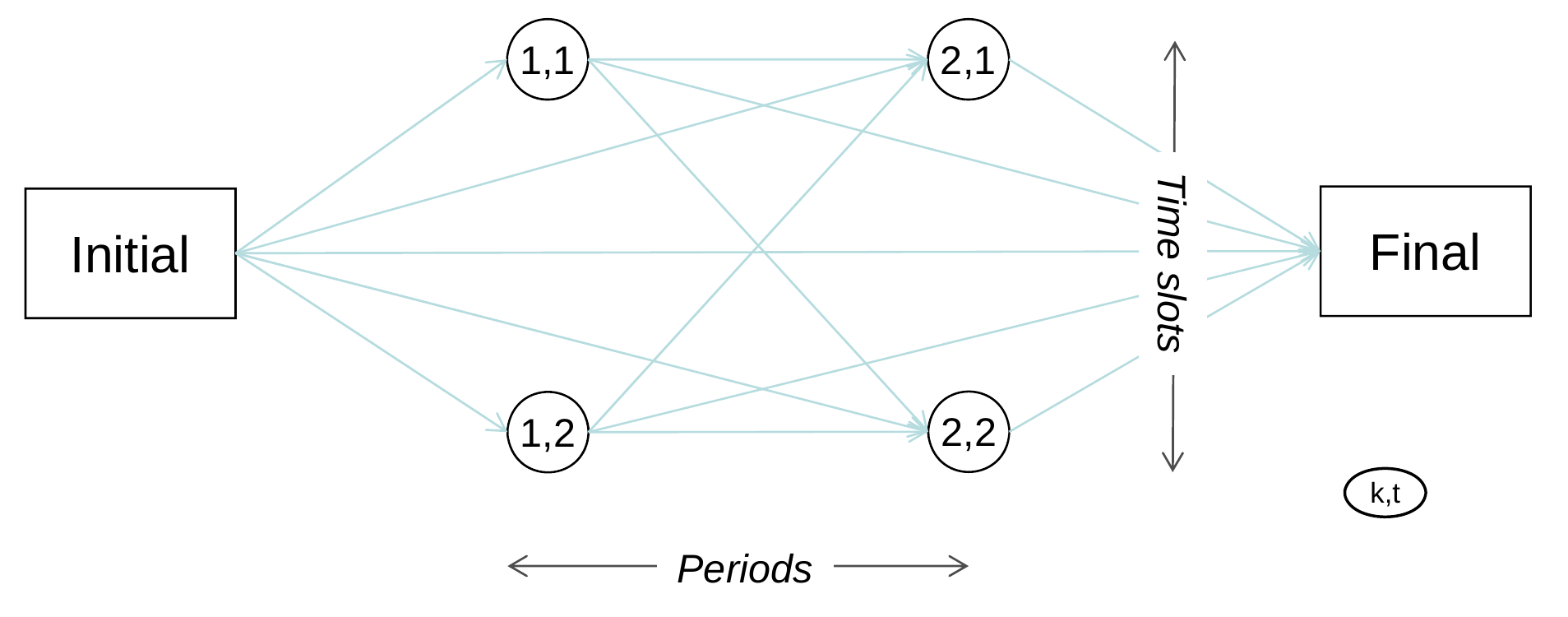}
\caption{\small Representation of a patient appointment scheduling
environment with two periods and two time slot per period. \label{f:net_flow}}
\end{figure}

The appointment scheduling environment can be represented by means
of a layered network. Each layer represents a period and each node
within a layer represents a time slot, i.e., node $(k,t)$ represents
time slot $t$ in period $k$.  A layered network with two periods and
two time slots per period is depicted in Figure \ref{f:net_flow}.
An arc $((k_1, t_1), (k_2, t_2))$ in the layered network represents
the option to schedule an appointment for a patient in time slot
$t_1$ in period $k_1$ followed by an appointment in time slot $t_2$
in period $k_2$. We assume that each patient is seen in period $0$.
While this assumption simplifies the modeling, it can be relaxed
easily.

A patient appointment schedule, i.e., the periods and the time slots within these periods in which the patient is scheduled to be seen, can be represented as a path in the layered network. The appointment scheduling problem is to determine patient appointment schedules
that maximize the aggregate health status of all patients over the planning horizon and in
which no time slot is assigned to more than one patient. A population appointment schedule, i.e., the set of schedules for all patients in the population, can be represented as a set of node disjoint paths in the layered network.

\subsection{Patient disease progression} \label{sb:dis_prog}

The health state of an asthmatic patient is defined by two factors:
severity and control. Severity can be interpreted as the intrinsic
susceptibility of the patient (a factor measured at a patient's
first appointment and a factor that does not change over time).
Severity creates different classifications of patients; the common
severity levels are mild intermittent, mild persistent, moderate
persistent, and severe persistent; see \citet{NHLBI2007guidelines}. Control is the extent to which a patient's asthma is under control,
and may change over time with treatment and natural disease
progression.  Categories for control vary within the asthma
community; however, MCF uses one category for controlled and three
sub-categories for uncontrolled, depending on the degree to which
the patient's asthma is not controlled.

\cite{Deoetal2011} characterize disease progression by a patient's
severity, the control state diagnosed and treatment performed at the
last visit, and the time since the last visit.  The authors model
disease progression as a Markov process.  Based on data from MCF,
the authors calibrate a per-period transition matrix $\bp$ to model
natural disease progression between control states for patients and
a transition matrix $\bq$ to represent the treatment effect of a
scheduled visit in terms of changing a patient's control status.
Recognizing that treatment is most effective just after a visit and
that natural disease progression occurs in subsequent periods, the
transition matrix $\bpq$ is applied to a patient's diagnosed control
state following a visit, and matrix $\bp$ is applied in following
periods until the next scheduled visit.

In our work, we focus on the special case with only two control states (\textit{$0$ = controlled} and \textit{$1$ = uncontrolled}) and ``perfect repair''. With perfect repair, a patient returns to the controlled state after a treatment, regardless of the state diagnosed at the visit.  After treatment, disease progression continues with the matrix $\bp$.  In this setting, disease progression can be characterized by severity (which influences $\bp$ as described below) and the time since the last visit.  As in \cite{Deoetal2011}, we assume that a patient's control cannot improve through natural disease progression, i.e., an uncontrolled patient can only become controlled through a scheduled visit).  With these assumptions, the disease progression matrix $\bp$ (for periods in which no visit is scheduled for a patient) and treatment matrix $\bq$ (for periods in which a visit is scheduled for a patient) are as follows:
\begin{align*}
    \bp = \left[ \begin{array}{cc}
    \alpha &1-\alpha\\
    0& 1
    \end{array}\right]~~~
    \bq = \left[ \begin{array}{cc}
    1 & 0\\
    1& 0
    \end{array}\right],
\end{align*}
where the first row and first column correspond to being in a controlled state and the second row and the second column correspond to being in an uncontrolled state.

\begin{figure}[htbp]
\vspace{-1cm}
        \centering
        \includegraphics[width=0.60\textwidth]{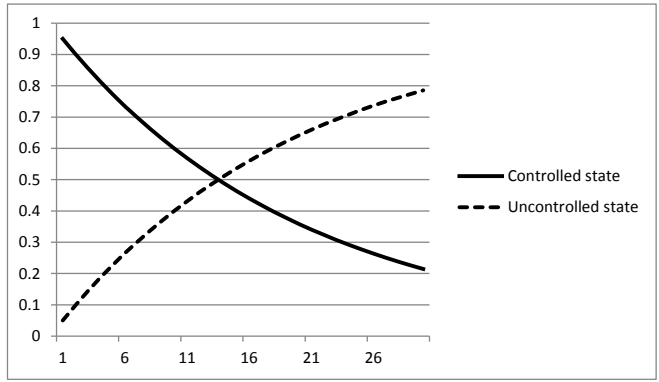}
\vspace{-1.5cm}\caption{\small Progression over time of the probability of
a patient's health state for $\alpha = 0.95$. \label{f:progression}}
\end{figure}

The parameter $\alpha$ represents the probability that a controlled
patient remains in a controlled health state in the following
period.  As shown in \cite{Deoetal2011}, the value of $\alpha$
depends on the patient's severity.  For ease of notation, we
formulate the optimization model using a single value of $\alpha$.
However, the computational study includes values that vary by
severity.  The probability that a controlled patient remains in the
controlled state decreases as the time since the last visit,
$\delta$, increases. More specifically, the probability that a
patient is in the controlled health state $\delta$ periods after his
last visit is $\alpha^\delta$, and, therefore, the probability that
a patient is in the uncontrolled health state $\delta$ periods after
his last visit is $1-\alpha^\delta$. (See Figure \ref{f:progression}
for an example with $\alpha = 0.95$.)

\subsection{Patient no-show probabilities}

Controlling the health states of patients is often complicated by
patients' lack of adherence to scheduled appointments. In the
context of MCF, this is due mostly to parents not showing
up at their child's appointment, which means the examination and
treatment of the child cannot occur, because a parent must be
present. Thus, associated with each patient $i \in \{1, ..., P\}$
and each time slot $t \in \{1, ..., T\}$, there is a patient no-show
probability $n_{it}$. These probabilities differ by time slot $t$
due to the relative convenience of the time slots (e.g., first
appointment of the day and during lunch time).  We assume that the
probabilities are the same in each period of the planning horizon,
although the model can be generalized to relax this assumption.

Determining patient no-show probabilities is challenging. As patients are seen infrequently, it is impossible to collect sufficient data to employ statistical techniques to estimate patient no-show  probabilities. However, a process can be put in place to get meaningful information from the patients themselves. At an intake consultation, initial time-of-day preference information has to be collected, and follow-up phone calls or emails have to take place at regular intervals to find out if the information on file is still accurate or whether time-of-day preferences have changed. Furthermore, if a no-show occurs, it is essential to assess whether inaccurate or out-of-date patient time-of-day preference information was a contributor, and, if so, take the necessary corrective actions.


\subsection{Patient appointment schedule evaluation}

In Section \ref{sb:pat_app_sch}, we discuss how a patient
appointment schedule can be represented as a path in a layered
network, where an arc in the path links two consecutive scheduled
appointments for a patient, and, in Section \ref{sb:dis_prog}, we
discuss how asthma control deteriorates with the time between
visits (the length of an arc). In this section, we present an
approach to evaluate patient schedules based on the probability of disease control over the planning period (to be
defined precisely next). Using disease control as an indicator of the quality of an appointment schedule
seems appropriate, as \cite{briggs2006cost}, for example, link
asthma control level to a health related quality of life and
\cite{price2002development} demonstrate the links between asthma
control and attack occurrence.

First, we consider the situation with perfect schedule
adherence, i.e., $n_{it} = 0$ for all $i \in \{1, ..., P\}$ and $t
\in \{1, ...,T\}$. In this setting, we associate the quantity
\begin{equation}
\sum_{\delta=1}^{\Delta} (1-\alpha^\delta),
\end{equation}
with an arc between two consecutive appointments that are $\Delta$ periods apart. We refer to this quantity, which is the sum of probabilities that a patient is in the uncontrolled state in the periods between the two appointments, as the \textit{aggregate probability} (realizing that the value, in fact, does not represent an actual probability). The aggregate probability $U_i$ of patient
$i$ being in an uncontrolled state during the planning horizon is
then simply the sum of the aggregate probabilities associated with
the arcs in the path representing the patient's appointment
schedule.

However, when patient schedule adherence is not perfect, i.e.,
$n_{it} > 0$ for some or all $i \in \{1, ..., P\}$ and $t \in \{1,
...,T\}$), the aggregate probability associated with an arc in the
path can no longer be calculated without knowledge of prior
appointments. The time between consecutive visits of a patient is no
longer equal to the number of periods represented by the length of
an arc, since there is a positive probability that the patient did
not show up at the appointment at the tail of the arc. The presence
of no-show probabilities implies that the time between visits is
uncertain, and thus calculating the aggregate probability that a
patient is in the uncontrolled state during the planning horizon
becomes more involved.

\begin{figure}[htbp]
    \centering
        \includegraphics[width=0.65\textwidth]{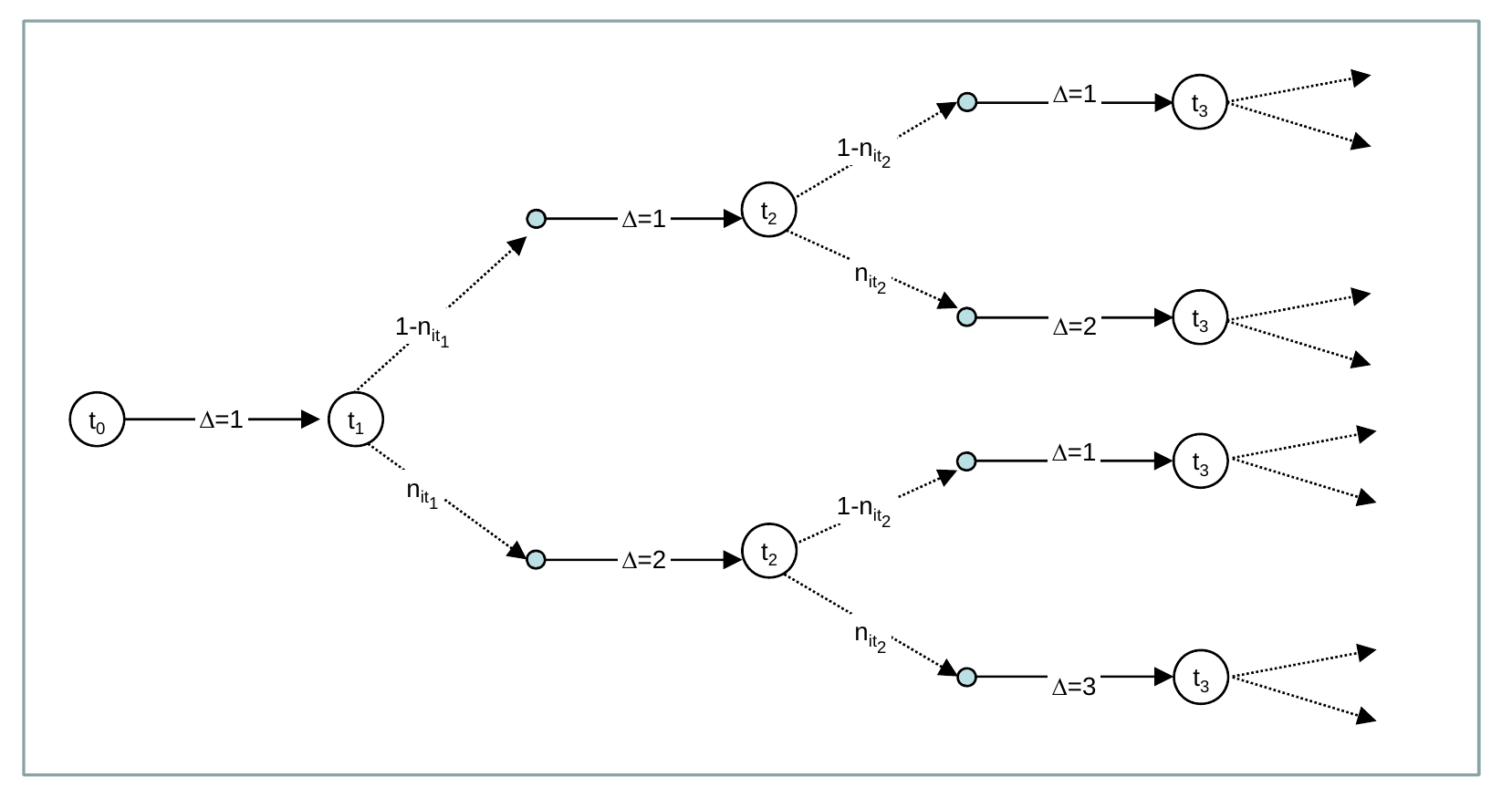}
    \caption{\small Time-since-last-visit probability tree for a patient appointment schedule $t_1$, $t_2$, $t_3$, $\dots$. \label{f:decision_tree}}
\end{figure}

To simplify the calculations, we model the option of not scheduling patient $i$ in a given period with a fictitious time slot $T+1$ with no-show probability $n_{i,T+1} = 1$, i.e., a patient that is scheduled \textit{not} to be seen in a period, will not be seen in that period with probability 1. By adding an additional node corresponding to this artificial time slot to the layered network (as well as the necessary arcs), a patient appointment schedule is represented by a path of exactly $K+1$ arcs, each connecting one period to the next. (Note that we allow more than one patient to be scheduled in this fictitious time slot.) We can now construct a time-since-last-visit probability tree for a patient appointment schedule. Figure \ref{f:decision_tree} presents the early periods of a
time-since-last-visit probability tree for a patient appointment
schedule with appointment time slots $t_1$, $t_2$, $t_3$, $\dots$.
As shown in Figure \ref{f:decision_tree}, the probability of the
time since the last visit, $l$, can be calculated explicitly at each
period (where, for presentational convenience, we have indicated the
time, $\Delta$, since the last visit on the arc into a node). Given
the assumption that all patients are seen in period 0, there are
$K+1$ possible values for the time since the last visit $l$. To
calculate the expected time since the last visit $l$, the
probability distribution function of all possible values of $l$ is
needed. Let $P_{i}^{kl}$ denote the probability that the number of
periods since the last visit of patient $i$ is $l$ immediately after
the scheduled appointment in period $k$, i.e.,
\begin{equation} \label{e:state_prob}
P_{i}^{kl} = \begin{cases}
   1 - n_{it} &\text{ if } l = 0\\
   n_{it} P_i^{k-1,l-1}   &\text{ otherwise, }
    \end{cases}
\end{equation}
assuming that the appointment in period $k$ is in time slot $t$. Since we assume that at the start of the planning horizon a patient has just been seen, we have $P_i^{00} = 1$ and  $P_i^{0l} = 0, \forall l > 0$. Recall that if patient $i$ is not scheduled in
period $k$, then $n_{it} = n_{i,T+1} = 1$.  Since the no-show
probability of the time slot impacts the control state of a patient
in the interval following the scheduled appointment, the control
state is calculated including the interval following the scheduled
appointment. With the distribution function defined at each
period $k$, the expected aggregate probability $E[U_i^k]$ of patient
$i$ being in an uncontrolled state after a visit in period $k$, including the interval following period $k$, is
\begin{equation} \label{e:exp_cost}
E[U_i^k] = E[U_i^{k-1}] + \sum_{l=0}^{k}  P_i^{kl} (1-\alpha^{l+1}),
\end{equation}
where $E[U_i^0] = 1 - \alpha$.

\section{Patient appointment scheduling approaches} \label{s:sch_app}

As mentioned in the introduction, the no-show rates at MCF vary by time of day, which suggests that taking time-of-day preferences into account during the construction of appointment schedules may be beneficial and may improve population health outcomes. As a consequence, the central questions underlying our research are whether time-of-day preference information can be incorporated in patient scheduling algorithms and whether the benefits of employing such algorithms can be quantified. To answer these questions, we develop and deploy a sophisticated appointment scheduling method, which considers individual patient characteristics (Section \ref{s:sol_app}), as well as relatively
simple and easy-to-use appointment scheduling methods, which divide the patients into groups with similar characteristics and schedule the patients within a group using a round-robin scheme, which we refer to as ``cohort scheduling policies'' (Section
\ref{s:rule-methods}).

\subsection{Optimization-based appointment scheduling methods} \label{s:sol_app}

\subsubsection{Model formulation}

Let the set of all possible patient appointment schedules be denoted
by $\mathcal{R}$.  Furthermore, let $b_{kt}^r$ for 
$k \in \{1,...,K\}$, $t \in \{1,...,T\}$, and $r \in \mathcal{R}$ indicate
whether or not a patient is seen in time slot $t$ in period $k$ in schedule
$r$ ($b_{kt}^r = 1$) or not ($b_{kt}^r = 0$), and let $u_i^r$ for $i
\in \{1,...,P\}$ and $r \in \mathcal{R}$ denote the expected
aggregate probability  $E[U_i^K]$ of being in an uncontrolled state
over the planning horizon when schedule $r$ is assigned to patient
$i$ (calculated using (\ref{e:exp_cost})). Finally, let $x_i^r$ for
$i \in \{1,...,P\}$ and $r \in \mathcal{R}$ be a binary variable
representing whether or not schedule $r$ is assigned to patient $i$
($x_i^r = 1$) or not ($x_i^r = 0$). Recall that we model the option
of not scheduling a patient in a given period with a time slot $T+1$
with capacity $C_{T+1} = P$; the capacity of all other time slots is
1. The optimization model is defined as
\begin{subequations}
\label{e:set_part}
\begin{equation}                \label{e:set_part_obj}
\min \sum_{r \in \mathcal{R}} \sum_{i = 1}^{P} u_i^r x_i^r
\end{equation} \indent
subject to
\begin{align}
\sum_{r \in \mathcal{R}} x_i^r &= 1     & i \in \{1, ..., P\}                \label{e:patient_assign} \\
\sum_{i = 1}^{P} \sum_{r \in \mathcal{R}} b_{kt}^r x_i^r & \leq C_t  & k \in \{1, ..., K\}, t \in \{1, ..., T+1\}          \label{e:slot_assign} \\
x_i^r &\in \{0,1\}                           & r \in \mathcal{R}, i
\in \{1, ..., P\}. \label{e:integer}
\end{align}
\end{subequations}

Rather than enumerating all possible patient appointment schedules
upfront, we use column generation to solve the linear programming
relaxation of (\ref{e:set_part}) and iteratively add new appointment
schedules to a restricted master problem \citep{BJNSV1998, dds2005}.
We relax constraints (\ref{e:patient_assign}) to $\sum_{r \in
\mathcal{R}} a_i^r x_i^r \geq 1$ for computational efficiency. Since
all patient appointment schedules have a positive aggregate
probability of being in an uncontrolled state, this will not change
the optimal solution.

We initialize the restricted master problem with the patient
schedules derived from a simple rotation policy (see Section
\ref{s:rule-methods}). After solving the linear programming
relaxation of (\ref{e:set_part}), we use a branch-and-bound approach
to obtain an integer solution.  We do not generate additional
columns throughout the branch-and-bound tree.

\subsubsection{Pricing problem formulation}

Given an optimal solution to the linear programming relaxation of
the restricted master problem, a pricing problem is solved to
determine whether there are any patient appointment schedules with
negative reduced costs. This can be done independently for each
patient.

\begin{figure}[htbp]
    \centering
        \includegraphics[width=0.9\textwidth]{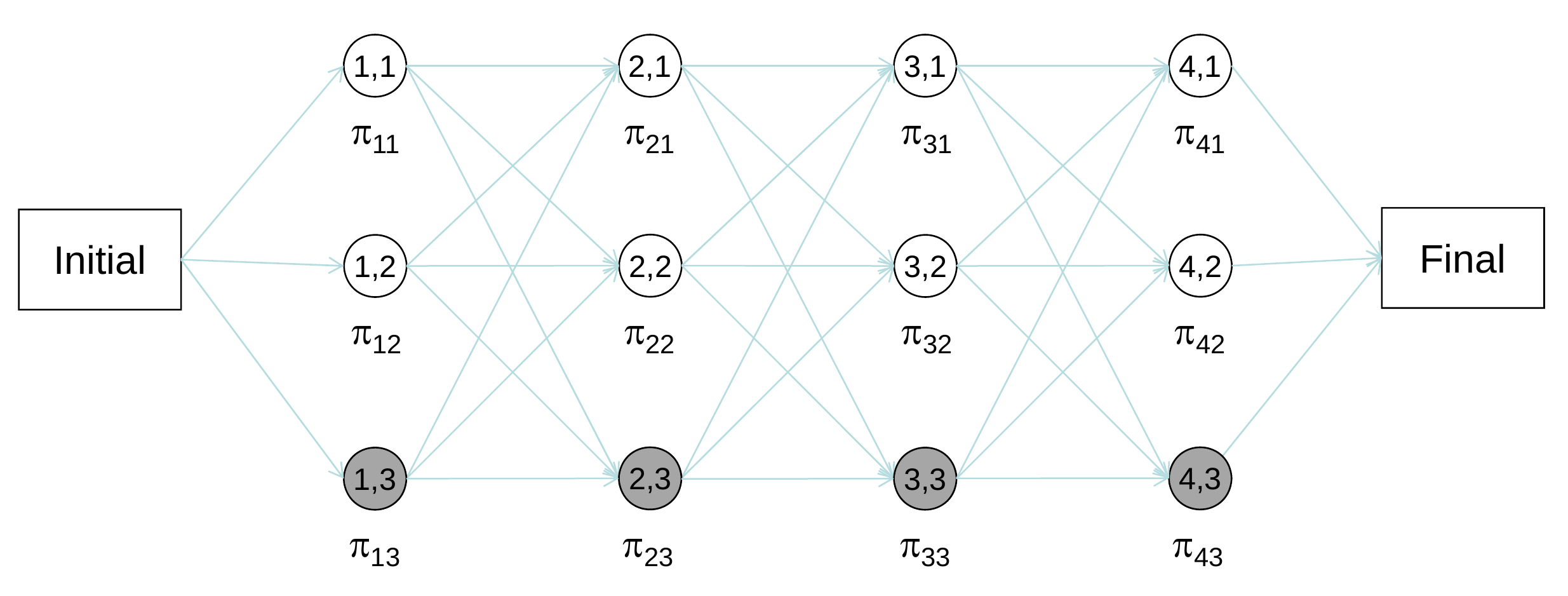}
    \caption{\small A layered network for a pricing problem for a single patient; 4 periods and 2 time slots. \label{f:pricing_prob}}
\end{figure}

Recall that a patient appointment schedule can be represented as a
path in a layered network. Figure \ref{f:pricing_prob} shows an
example of the layered network for an environment with four periods
and two time slots per period. Note that each layer, corresponding
to a period, includes an additional node to account for the patient
not being scheduled in that period.

Let $\sigma_i$ denote the dual variable associated with the
relaxation of constraint (\ref{e:patient_assign}) for patient $i$
and let $\pi_{kt}$ denote the dual variable associated with
constraint (\ref{e:slot_assign}) for  period $k$ and time slot $t$.
The reduced cost of an appointment schedule for patient $i$ is given
by the expected aggregate probability of being in an uncontrolled
state of that appointment schedule plus the sum of the dual values
associated with the time slots in that appointment schedule and the
dual value associated with the constraint that ensures exactly one
appointment schedule is selected for the patient.  This is
equivalent to the value of the corresponding path in the layered
network plus the sum of the dual values associated with the nodes
visited on that path and the dual value associated with the
constraint that ensures exactly one appointment schedule is selected
for the patient (this last term is independent of the path in the
layered network). Therefore, determining whether a patient
appointment schedule with negative reduced cost exists for patient
$i$ can be done by solving a shortest path problem on the layered
network.

The {\it adjusted} expected aggregate probability $\hat{E}[U_i^k]$
of a partial path for patient $i$ ending in node $(k,t)$, i.e.,
adjusted by the dual values associated with the nodes visited on the
partial path, is given by
\begin{equation} \label{e:exp_red_cost}
\hat{E}[U_i^k] = \hat{E}[U_i^{k-1}] + \sum_{l=0}^{k} P^{k,l}_i
(1-\alpha^{l+1}) + \pi_{kt},
\end{equation}
which involves the (discrete) probability distribution of the time
since the last visit, represented by the $K+1$ dimensional vector
$(P^{k,0}_i, P^{k,1}_i, \ldots, P^{k,K}_i)$ and defined by (\ref{e:state_prob}).
Note that in (\ref{e:exp_red_cost}) we use the fact that in
period $k$ the probability that the time since the last visit is
greater than $k$ is zero (because every patient is seen in period
0).

The value of the dual variable associated with the relaxation of
constraints (\ref{e:patient_assign}) that ensure exactly one
appointment schedule is selected for patient $i$, i.e., $\sigma_i$,
is added to the adjusted expected aggregate probability of a
(complete) path to determine the reduced cost of the path. The
pricing problem finds for each patient $i \in \{1, ..., P\}$ a path
with minimum reduced cost and adds the corresponding column to the
restricted master problem if the reduced cost of that path is
negative. The restricted master problem is resolved to obtain a new
optimal dual solution and the process repeats as long as any columns
with negative reduced costs are found.

\subsubsection{Pricing problem solution approaches}

With the inclusion of the ``no appointment'' node, the network has a
simple layered structure in which each layer corresponds to a
period and in which there are only arcs between consecutive layers.
Thus, any path from the source to the sink visits exactly one node
in each layer; see Figure \ref{f:pricing_prob} for an example. The
structure of the layered network is the same for all patients; only
the no-show rates and the severity differ by patient. The dual
values change each time the pricing problem is solved.

Because solving the pricing problem optimally involves solving a
multi-label shortest path problem with $K+2$ labels ($K+1$ for the
probability vector and one for the adjusted expected cost), solving
the pricing problem for large values of the planning horizon $K$ can
become prohibitive. Therefore, we consider the following heuristic
for solving the pricing problem. Rather than using the probability
distribution of the time since the last visit, we use the expected
time $E[\Delta^k]$ since the last visit, which can be calculated as
follows: $E[\Delta^k] = (1-n_{it})(1) + n_{it} (E[\Delta^{k-1}] +
1)$. This reduces the number of labels to maintain in the
multi-label shortest path problem from $K+2$ to $2$.  Computational
tests show that the heuristic produces near-optimal solutions in an
acceptable amount of time.

\subsection{Cohort-based appointment scheduling methods}
\label{s:rule-methods}

In this section, we present cohort-based scheduling methods, which (1)
partition the patient population into cohorts based on some 
differentiating factor, e.g., time-of-day preference, disease
severity, and reliability; (2) use a simple rule to assign time slots in
the planning period to each of the cohorts; and (3) apply a
simple rotation policy to assign time slots to the patients in a
cohort. Cohort-based scheduling methods are intuitive and
easy-to-use and are quite effective when a heterogeneous patient
population can easily be partitioned into cohorts.

\subsubsection{Cohort development}  \label{sbb:cohortdevel}

In this study, we consider three differentiating factors for
grouping patients into cohorts: time-of-day preference, disease
severity, and reliability. Cohort strategies can be characterized by
the number of factors considered for differentiating patients and
the specific features used. A 1-level cohort strategy based on
time-of-day preference, for example, partitions the patient
population into two or more cohorts based on patients' time-of-day
preferences, e.g., a cohort that prefers morning time slots, a
cohort that prefers noon-time time slots, and a cohort that prefers
afternoon time slots. The 0-level cohort strategy has a single cohort consisting of the entire patient population and thus does not distinguish patients and treats all patients the same.

\subsubsection{Allocating time slots to cohorts} \label{sbb:timeassign}

Due to the natural relation between time-of-day preferences and time slots, the logic for allocating time slots to cohorts determined using time-of-day preferences should be different from the logic for allocating time slots to cohorts determined using either disease severity or reliability.

First, we consider cohorts determined using time-of-day preferences. When patients are partitioned into a morning cohort and an afternoon cohort, the
number of patients in each cohort is roughly equal, and there are eight time slots in a period (as is the case in our computational experiments), then assigning the 4 morning time slots to the morning cohort and the 4 afternoon
time slots to the afternoon cohort is the natural course of action.
On the other hand, if there are three times as many patients in the
afternoon cohort, then assigning the first 2 time slots (morning) to
the morning cohort and the last 6 time slots (late morning and
afternoon) to the afternoon cohort is the natural course of action.
In more complex settings, where the numbers do not work out as
nicely as in the examples above, more sophisticated approaches can
be employed, somewhat similar to the one discussed next for allocating slots
to cohorts created by differentiating based on disease severity and
reliability.

For ease of presentation, we assume that there are two cohorts with $n_1$ and $n_2$ patients ($n_1 \leq n_2$), respectively, and that a total of $m = |K||T|$ time slots have to be assigned to the patients over the planning period ($m \gg n_1 + n_2$). The first step is to divide the $m$ time slots over the two cohorts. One possibility is a proportional allocation according to the number of patients in the cohort, but in many cases this not the best choice. For example, when the cohorts are created based on disease
severity and the number of patients in the cohorts is the same, it
is probably better to allocate more time slots to the cohort with patients
with a higher severity level. For now, suppose that a fraction $f < 0.5$ of
the time slots is allocated to the first cohort, i.e., $\lceil f m
\rceil$ time slots will be used for appointments of patients in the
first cohort. We spread these time slots equally spaced over
the total $m$ time slots by allocating time slots $\lceil \frac{j}{f} \rceil$ for $j = 1, \ldots, \lceil f m \rceil$ to the first cohort. The
remaining time slots are allocated to the second cohort. For
example, if there are 80 time slots in the planning period and 40\% of them
($f = 0.4$) are allocated to the first cohort, then the first cohort
will get time slots $3 = \lceil \frac{1}{0.4} \rceil$, $5 = \lceil
\frac{2}{0.4} \rceil$, $8 = \lceil \frac{3}{0.4} \rceil$, etc. The
scheme can easily be extended to accommodate more than two cohorts
by applying the above procedure recursively, e.g., allocate time
slots to the cohort with the smallest fraction of time slots, then
allocated time slots to the cohort with the second smallest fraction 
of time slots, etc. Thus, to define a cohort strategy, one only needs to
decide on the fraction of slots that will be allocated to each of
the cohorts. Once that decision is made the time slots are
allocated automatically.

\subsubsection{Scheduling patients within each cohort}  \label{sbb:rotation}

We use a simple rotation policy to assign patients in a cohort to
the time slots allocated to the cohort, i.e., a round-robin
scheduling rule which schedules patients by patient index. See
Algorithm \ref{alg:rot} for a more precise description.
\begin{algorithm}
	\DontPrintSemicolon
    $i \leftarrow 1$ \;
    \For{$k\leftarrow 1$ \KwTo $K$}{
	    \For{$t\leftarrow 1$ \KwTo $T$}{
	        Assign patient $i$ to time slot $t$ in period $k$ \;
            \If{$i = P$}{$i \leftarrow 1$} \Else{$i \leftarrow i+1$}
        }
	}
\caption{\small Creating an appointment schedule with a rotation policy. \label{alg:rot}}
\end{algorithm}

The rotation policy has the advantage that it automatically spreads
out appointments and diversifies the time slots of the appointments
(unless the number of time slots in a period is a divider of the
number of patients in a cohort).  When the number of patients is a
multiple of the number of time slots, a diversification can be
introduced by using a slot-reversing rotation policy; see the
Appendix A for details.

\section{Computational study} \label{s:comp_study}

We have conducted an extensive computational study to (1) assess the benefit of considering time-of-day preferences when scheduling appointments (by incorporating no-show probabilities during schedule creation), and (2) assess the qualitative differences between the optimization-based appointment scheduling method and
the simpler and easier-to-use cohort-based appointment scheduling
methods.

\subsection{Instances}

To assess the benefit of considering patient time-of-day preferences during appointment scheduling and, more generally, accounting for different patient characteristics during appointment scheduling, we create a set of instances with varying patient profiles along the key dimensions of severity, reliability, and time-of-day preferences. Each instance has 20 patients, covers a planning horizon of 13 periods, and each period has 8 time slots.

\subsubsection*{Time-of-day preferences}

Time-of-day preferences are modeled in terms of no-show probabilities (e.g., low no-show probabilities for morning time slots indicate a preference for morning time slots). We consider three categories of time-of-day preference: \textit{AM}, \textit{Noon}, and \textit{PM}.  The no-show probabilities associated with each of these time-of-day preferences, for both a strong preference variant and a weak preference variant, are shown in Table \ref{t:no_show}.  We consider six patient population profiles, shown in the right-most part of Table \ref{t:no_show}, in which the number of patients with a specific time-of-day preference differs; Profile I \& II: homogeneous or almost homogeneous preferences; Profile III \& IV \& V: mixed AM, Noon, and PM preferences; Profile VI: balanced AM, Noon, and PM preferences.

\begin{table}[htbp]
\caption{\small Time-of-day preferences and slot-dependent no-show probabilities.}
\vspace{6pt} \label{t:no_show} \centering
\scriptsize
\begin{tabular}{l|l| r r r r r r r r | r r r r r r}
\toprule
Preference & Category & \multicolumn{8}{c}{No-show probability}  & \multicolumn{6}{|c}{Profiles} \\
Strength & & \multicolumn{8}{c}{(time slot)} & \multicolumn{6}{|c}{(\# patients)} \\
\hline
& & 1 & 2 & 3 & 4 & 5 & 6 & 7 & 8 & I & II & III & IV & V & VI \\
\hline
Strong & \textit{AM}      &   0.05  &   0.05  &   0.05   &   0.05   &   0.35  &   0.35  &  0.35  &   0.35   & 20 & 16 & 10 & 10 & 5  & 7\\
& \textit{Noon}    &   0.35  &   0.35  &   0.35   &   0.05   &   0.05  &   0.35  &  0.35  &   0.35   & 0  &  2 &  5 &  0 & 10 & 6\\
& \textit{PM}      &   0.35  &   0.35  &   0.35   &   0.35   &   0.05  &   0.05  &  0.05  &   0.05   & 0  &  2 &  5 & 10 & 5  & 7\\
\hline
Weak & \textit{AM}      &   0.05  &   0.05  &   0.05   &   0.05   &   0.15  &   0.15  &  0.15  &   0.15   & 20 & 16 & 10 & 10 & 5  & 7\\
& \textit{Noon}    &   0.15  &   0.15  &   0.15   &   0.05   &   0.05  &   0.15  &  0.15  &   0.15   & 0  &  2 &  5 &  0 & 10 & 6\\
& \textit{PM}      &   0.15  &   0.15  &   0.15   &   0.15   &   0.05  &   0.05  &  0.05  &   0.05   & 0  &  2 &  5 & 10 & 5  & 7\\
\bottomrule
\end{tabular}
\normalsize
\end{table}%

\subsubsection*{Severity}

Patient severity is modeled with different values of $\alpha$, the probability that a controlled patient remains in a controlled health state in the period following treatment. We consider severities in the range from \textit{Mild}, modeled with $\alpha = 0.9$, to \textit{Severe}, modeled with $\alpha = 0.8$. We consider four patient population profiles, shown in the right-most part of Table \ref{t:severity}, in which the number of patients with a specific severity level differs; Profile I: homogeneous mild; Profile II: mixed mild/severe; Profile III: homogeneous severe; Profile IV: varied with severities in the interval [0.8 - 0.9] (i.e., between mild and severe). In Profile II, the patients in the population with a mild severity level are selected randomly (and thus so are the patients with a severe severity level). As a consequence, the number of mild and severe patients with a similar time-of-day preferences might not be balanced, e.g., if 10 patients have a morning time-of-day preference, it is possible that three have a mild level of severity and seven have a severe level of severity. The severity of the patients in the varied profile is drawn randomly from a uniform distribution with lower and upper bounds 0.8 and 0.9, respectively.
\begin{table}[htbp]
\caption{\small Severity profiles.} \vspace{4pt}
\label{t:severity} \centering
\begin{tabular}{l | c | r r r r  }
	\toprule
Category              & \multicolumn{1}{c}{Control probability}  & \multicolumn{4}{|c}{Profiles} \\
                      & \multicolumn{1}{c}{($\alpha$)}           & \multicolumn{4}{|c}{(\# patients)} \\

\hline
        &  & I & II & III & IV \\
\hline
Mild    &   0.8         & 20 & 10 &  0 & 0 \\
Varying &   [0.8 - 0.9] & 0  &  0 &  0 & 20 \\
Severe  &   0.9         & 0  & 10 & 20 & 0 \\
	\bottomrule
\end{tabular}%
\end{table}%

\subsubsection*{Reliability}

Patient reliability is modeled by adjusting the no-show probabilities associated with time-of-day preferences. We consider reliabilities in the range from \textit{Reliable}, modeled by multiplying the no-show probabilities associated with time-of-day preferences by 0.8, to \textit{Unreliable}, modeled by keeping the no-show probabilities associated with time-of-day preferences unchanged. We consider four patient population profiles, shown in the right-most part of Table \ref{t:reliability}, in which the number of patients with a specific reliability level differs; Profile I: homogeneous reliable; Profile II: mixed reliable/ unreliable; Profile III: homogeneous unreliable; Profile IV: varied with reliabilities in the interval [0.8 - 1.0] (drawn randomly from a uniform distribution with lower and upper bounds 0.8 and 1.0, respectively). Again, reliable and unreliable patients in the \textit{mixed} profile are selected randomly, and thus the number of reliable and unreliable patients with a similar time-of-day preference (and severity level) might not be balanced.
\begin{table}[htbp]
	\caption{\small Reliability profiles.} \vspace{4pt}
	\label{t:reliability} \centering
	\begin{tabular}{l | c | r r r r  }
		\toprule
		Category              & \multicolumn{1}{c}{Reliability}  & \multicolumn{4}{|c}{Profiles} \\
		& & \multicolumn{4}{|c}{(\# patients)} \\
		
		\hline
		&  & I & II & III & IV \\
		\hline
		More reliable   &   0.8         & 20 & 10 &  0 & 0 \\
		Varying    &   [0.8 - 1.0] & 0  &  0 &  0 & 20 \\
		Less reliable &   1.0         & 0  & 10 & 20 & 0 \\
		\bottomrule
	\end{tabular}%
\end{table}%

\vspace{12pt} \noindent We have combined the above patient population profiles into a total of 42 instances as shown in the first five columns of  Table \ref{t:cohort_perf}. For each of these instances, we examine two variants, one in which the time-of-day preferences are strong, and one in which the time-of-day preferences are weak.

\subsection{Computational results}

To evaluate the benefit of accounting for different patient characteristics during appointment scheduling, and to determine the effectiveness of cohort policies, we evaluate the expected aggregate control probability $z$ for the patient population of the appointment schedule produced by a particular cohort strategy, relative to the expected aggregate control probability $z^*$ for the patient population of the appointment
schedule produced by the optimization-based method.

In Table \ref{t:cohort_perf}, we present the percentage performance gap, $100\frac{z-z^*}{z^*}$. The first five columns of Table \ref{t:cohort_perf} describe the instance characteristics.  Columns 6-9 present the performance gaps for instances in which patients have strong time-of-day preferences and columns 10-14 present performance gaps for instances in which patients have weak time-of-day preferences.  Recall that the 0-level cohort policy is a simple rotation policy that does not distinguish patients and treats all patients the same. The three 1-level cohort scheduling policies relate to the distinguishing factor used to define the cohorts, i.e., time-of-day preference (T), disease severity (S), and reliability (R). A 1-level cohort policy requires the specification of the number of slots allocated to each of the cohorts. These allocations can be found in Tables \ref{coh:tod}, \ref{coh:sev}, and \ref{coh:rel} in Appendix B, respectively. We note that because the optimization-based method uses heuristic pricing and does not generate additional columns during the tree search, it is possible to see negative gaps (indicating that a cohort-strategy has produced a better solution).

\begin{table}[htbp]
  \centering
  \caption{\small Performance of cohort scheduling policies relative to the optimization-based approach \label{t:cohort_perf}} \vspace{6pt}
	\scriptsize
	\centering
    \begin{tabular}{r|rr|rr|r|rrr|r|rrr}
    \toprule
    \multicolumn{5}{c|}{Instance characteristics}      & \multicolumn{4}{c|}{Strong time preference (0.05 v 0.35)} & \multicolumn{4}{c}{Weak time preference (0.05 v 0.15)} \\
    \midrule
    \multicolumn{1}{c|}{ToD Preference} & \multicolumn{2}{c|}{Severity} & \multicolumn{2}{c|}{Reliability} & 0-level & \multicolumn{3}{c|}{1-level} & 0-level & \multicolumn{3}{c}{1-level} \\
    \multicolumn{1}{c|}{Profile} & \multicolumn{1}{c}{0.8} & \multicolumn{1}{c|}{0.9} & 0.8   & 1     &       & \multicolumn{1}{c}{T} & \multicolumn{1}{c}{S} & \multicolumn{1}{c|}{R} &       & \multicolumn{1}{c}{T} & \multicolumn{1}{c}{S} & \multicolumn{1}{c}{R} \\ \hline
    \multicolumn{1}{c|}{\multirow{3}[2]{*}{I}} & 20    & 0     & \multicolumn{2}{c|}{\multirow{3}[1]{*}{Random}} & 1.6\% &       &       & 1.9\% & 0.7\% &       &       & 0.9\% \\
    \multicolumn{1}{c|}{} & 10    & 10    & \multicolumn{2}{c|}{} & 1.3\% &       & 0.6\% & 1.4\% & 2.1\% &       & 0.9\% & 2.4\% \\
    \multicolumn{1}{c|}{} & 0     & 20    & \multicolumn{2}{c|}{} & 2.0\% &       &       & 2.5\% & 1.0\% &       &       & 1.3\% \\ \cline{2-13}
    \multicolumn{1}{c|}{AM:20}  & \multicolumn{2}{c|}{\multirow{4}[1]{*}{Random}} & 0     & 20    & 1.1\% &       & 2.7\% &       & 1.4\% &       & 1.8\% &  \\
    \multicolumn{1}{c|}{Noon:0} & \multicolumn{2}{c|}{} & 10    & 10    & 2.4\% &       & 3.8\% & 4.7\% & 1.5\% &       & 1.9\% & 3.3\% \\
    \multicolumn{1}{c|}{PM:0} & \multicolumn{2}{c|}{} & 20    & 0     & 1.1\% &       & 2.3\% &       & 1.3\% &       & 1.7\% &  \\
    \multicolumn{1}{c|}{} & \multicolumn{2}{c|}{} & \multicolumn{2}{c|}{Random} & 1.8\% &       & 3.3\% & 2.4\% & 1.5\% &       & 1.9\% & 1.9\% \\ \hline
    \multicolumn{1}{c|}{\multirow{3}[2]{*}{II}} & 20    & 0     & \multicolumn{2}{c|}{\multirow{3}[1]{*}{Random}} & 6.8\% & 0.9\% &       & 7.3\% & 2.3\% & 0.4\% &       & 2.6\% \\
    \multicolumn{1}{c|}{} & 10    & 10    & \multicolumn{2}{c|}{} & 6.8\% & 0.6\% & 5.3\% & 7.4\% & 3.5\% & 1.6\% & 1.9\% & 3.8\% \\
    \multicolumn{1}{c|}{} & 0     & 20    & \multicolumn{2}{c|}{} & 8.4\% & 1.2\% &       & 9.3\% & 2.9\% & 0.5\% &       & 3.3\% \\ \cline{2-13}
    \multicolumn{1}{c|}{AM:16} & \multicolumn{2}{c|}{\multirow{4}[1]{*}{Random}} & 0     & 20    & 8.0\% & 0.7\% & 8.9\% &       & 3.4\% & 1.0\% & 3.6\% &  \\
    \multicolumn{1}{c|}{Noon:2} & \multicolumn{2}{c|}{} & 10    & 10    & 8.2\% & 1.3\% & 8.9\% & 9.0\% & 3.2\% & 0.9\% & 3.4\% & 4.5\% \\
    \multicolumn{1}{c|}{PM:2} & \multicolumn{2}{c|}{} & 20    & 0     & 6.5\% & 0.6\% & 7.1\% &       & 2.9\% & 1.0\% & 3.1\% &  \\
    \multicolumn{1}{c|}{} & \multicolumn{2}{c|}{} & \multicolumn{2}{c|}{Random} & 7.6\% & 1.1\% & 8.4\% & 8.6\% & 3.3\% & 1.1\% & 3.5\% & 3.8\% \\ \hline
    \multicolumn{1}{c|}{\multirow{3}[2]{*}{III}} & 20    & 0     & \multicolumn{2}{c|}{\multirow{3}[1]{*}{Random}} & 11.0\% & 1.4\% &       & 11.9\% & 3.5\% & 0.4\% &       & 3.9\% \\
    \multicolumn{1}{c|}{} & 10    & 10    & \multicolumn{2}{c|}{} & 9.3\% & 2.4\% & 9.1\% & 9.4\% & 5.4\% & 2.0\% & 4.2\% & 5.7\% \\
    \multicolumn{1}{c|}{} & 0     & 20    & \multicolumn{2}{c|}{} & 13.9\% & 1.7\% &       & 15.2\% & 4.4\% & 0.5\% &       & 4.9\% \\ \cline{2-13}
    \multicolumn{1}{c|}{AM:10} & \multicolumn{2}{c|}{\multirow{4}[1]{*}{Random}} & 0     & 20    & 14.3\% & 2.0\% & 16.1\% &       & 5.1\% & 1.2\% & 5.7\% &  \\
    \multicolumn{1}{c|}{Noon:5} & \multicolumn{2}{c|}{} & 10    & 10    & 12.9\% & 1.7\% & 14.6\% & 16.3\% & 4.8\% & 1.2\% & 5.3\% & 6.9\% \\
    \multicolumn{1}{c|}{PM:5} & \multicolumn{2}{c|}{} & 20    & 0     & 11.4\% & 1.7\% & 12.8\% &       & 4.3\% & 1.2\% & 4.7\% &  \\
    \multicolumn{1}{c|}{} & \multicolumn{2}{c|}{} & \multicolumn{2}{c|}{Random } & 13.0\% & 2.1\% & 14.8\% & 14.2\% & 4.7\% & 1.2\% & 5.2\% & 5.3\% \\ \hline
    \multicolumn{1}{c|}{\multirow{3}[2]{*}{IV}} & 20    & 0     & \multicolumn{2}{c|}{\multirow{3}[1]{*}{Random}} & 10.2\% & -0.1\% &       & 10.3\% & 3.3\% & -0.1\% &       & 3.4\% \\
    \multicolumn{1}{c|}{} & 10    & 10    & \multicolumn{2}{c|}{} & 13.2\% & 2.1\% & 11.9\% & 13.5\% & 5.6\% & 2.1\% & 4.1\% & 5.8\% \\
    \multicolumn{1}{c|}{} & 0     & 20    & \multicolumn{2}{c|}{} & 12.7\% & -0.2\% &       & 12.9\% & 4.1\% & -0.2\% &       & 4.3\% \\ \cline{2-13}
    \multicolumn{1}{c|}{AM:10} & \multicolumn{2}{c|}{\multirow{4}[1]{*}{Random }} & 0     & 20    & 14.0\% & 1.0\% & 17.2\% &       & 5.3\% & 1.0\% & 6.2\% &  \\
    \multicolumn{1}{c|}{Noon:0} & \multicolumn{2}{c|}{} & 10    & 10    & 12.7\% & 1.0\% & 15.4\% & 14.7\% & 4.8\% & 1.0\% & 5.6\% & 6.5\% \\
    \multicolumn{1}{c|}{PM:10} & \multicolumn{2}{c|}{} & 20    & 0     & 11.4\% & 1.0\% & 13.8\% &       & 4.4\% & 1.0\% & 5.1\% &  \\
    \multicolumn{1}{c|}{} & \multicolumn{2}{c|}{} & \multicolumn{2}{c|}{Random } & 12.5\% & 1.0\% & 15.2\% & 12.9\% & 4.8\% & 1.0\% & 5.6\% & 5.2\% \\ \hline
    \multicolumn{1}{c|}{\multirow{3}[2]{*}{V}} & 20    & 0     & \multicolumn{2}{c|}{\multirow{3}[1]{*}{Random}} & 11.4\% & -0.1\% &       & 12.9\% & 3.7\% & 0.0\% &       & 4.2\% \\
    \multicolumn{1}{c|}{} & 10    & 10    & \multicolumn{2}{c|}{} & 14.7\% & 1.6\% & 16.6\% & 16.2\% & 6.0\% & 1.6\% & 5.4\% & 6.6\% \\
    \multicolumn{1}{c|}{} & 0     & 20    & \multicolumn{2}{c|}{} & 14.3\% & -0.1\% &       & 16.5\% & 4.6\% & -0.1\% &       & 5.3\% \\ \cline{2-13}
    \multicolumn{1}{c|}{AM:5} & \multicolumn{2}{c|}{\multirow{4}[1]{*}{Random }} & 0     & 20    & 15.6\% & 0.8\% & 18.2\% &       & 5.5\% & 0.8\% & 6.3\% &  \\
    \multicolumn{1}{c|}{Noon:10} & \multicolumn{2}{c|}{} & 10    & 10    & 13.8\% & 0.8\% & 16.2\% & 18.5\% & 5.0\% & 0.8\% & 5.7\% & 7.4\% \\
    \multicolumn{1}{c|}{PM:5} & \multicolumn{2}{c|}{} & 20    & 0     & 12.5\% & 0.8\% & 14.5\% &       & 4.6\% & 0.8\% & 5.2\% &  \\
    \multicolumn{1}{c|}{} & \multicolumn{2}{c|}{} & \multicolumn{2}{c|}{Random} & 13.7\% & 0.8\% & 16.3\% & 15.7\% & 4.9\% & 0.8\% & 5.7\% & 5.8\% \\ \hline
    \multicolumn{1}{c|}{\multirow{3}[2]{*}{VI}} & 20    & 0     & \multicolumn{2}{c|}{\multirow{3}[1]{*}{Random}} & 9.8\% & 1.0\% &       & 9.9\% & 3.1\% & 0.3\% &       & 3.3\% \\
    \multicolumn{1}{c|}{} & 10    & 10    & \multicolumn{2}{c|}{} & 12.2\% & 2.3\% & 13.0\% & 12.2\% & 5.3\% & 2.3\% & 4.4\% & 5.5\% \\
    \multicolumn{1}{c|}{} & 0     & 20    & \multicolumn{2}{c|}{} & 11.3\% & 0.1\% &       & 11.6\% & 3.6\% & 0.1\% &       & 3.9\% \\ \cline{2-13}
    \multicolumn{1}{c|}{AM:7} & \multicolumn{2}{c|}{\multirow{4}[1]{*}{Random }} & 0     & 20    & 12.1\% & 0.9\% & 14.2\% &       & 4.7\% & 1.2\% & 5.4\% &  \\
    \multicolumn{1}{c|}{Noon:6} & \multicolumn{2}{c|}{} & 10    & 10    & 11.6\% & 1.5\% & 13.5\% & 14.8\% & 4.5\% & 1.4\% & 5.1\% & 6.5\% \\
    \multicolumn{1}{c|}{PM:7} & \multicolumn{2}{c|}{} & 20    & 0     & 9.7\% & 1.0\% & 11.4\% &       & 4.0\% & 1.3\% & 4.5\% &  \\
    \multicolumn{1}{c|}{} & \multicolumn{2}{c|}{} & \multicolumn{2}{c|}{Random} & 11.3\% & 1.5\% & 13.4\% & 11.7\% & 4.4\% & 1.3\% & 5.1\% & 4.8\% \\
    \bottomrule
    \end{tabular}%
	\normalsize
\end{table}%

\vspace{12pt} \noindent \textbf{Analysis of the rotation policy
(0-level cohort scheduling policy)}

\begin{observation}
A simple rotation policy performs well only when all patients have the same time-of-day preferences.
\end{observation}

When patients have the same time-of-day preferences (i.e., time-of-day patient population profile I), the simple rotation policy results in a population appointment schedule with a level of aggregate control that is reasonably close to that of the population appointment schedule produced by the more sophisticated optimization-based approach (gaps of about one to two percent).  With the rotation policy, all patients are seen with the same frequency, and desirable and undesirable time slots are assigned alternatingly.

\begin{figure}[htbp]
	\centering
	\includegraphics[width=0.9\textwidth]{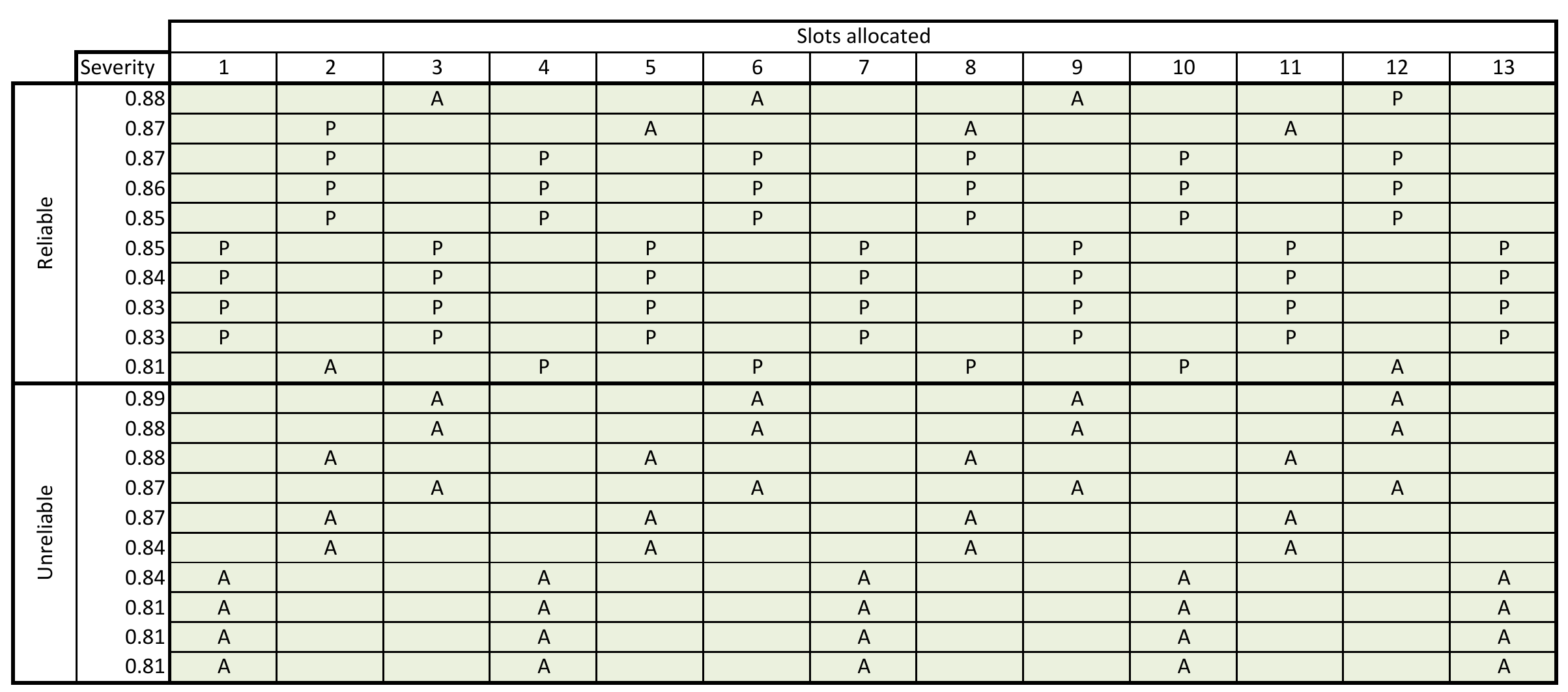}
	\caption{\small Time slots assigned in the optimization-based solution: Profile I: homogeneous AM time slot; mixed reliability; varying severity: (A) desirable morning slot (1-4); (P) less desirable afternoon slot (5-8). \label{f:I_MSR_VS}}
\end{figure}

The optimization-based approach exploits the full flexibility of assigning slots, i.e., the number of slots to assign to a patient, the specific periods in which to assign a slot to a patient (the spread), and the type of slot to assign to a patient (desirable or undesirable), and considers all patients in the population simultaneously. As a result, a (slightly) better population appointment schedule is obtained even in this setting.
Figure \ref{f:I_MSR_VS} shows the 13-period population appointment schedule produced by the optimization-based approach when all patients have a strong preference for morning slots, patients have mixed reliability, and severity levels vary across patients.
The first ten rows show the patient appointment schedules for the reliable patients and the second ten rows for the less reliable patients. Within each reliability group, patients are shown in nondecreasing order of severity. Note that this means, in some sense, that the patients that require the most carefully constructed appointment schedules appear in the bottom rows and the patients for which there is more leeway in constructing their appointment schedules appear in the top rows.

An examination of the population appointment schedule reveals the logic ``applied by'' the optimization-based approach. The patients that require the most carefully constructed appointment schedules (i.e., severe, but unreliable patients) are given their preferred slots and more of them when their disease is more severe. No afternoon slots are allocated to these patients. The patients for which there is more leeway in constructing their appointment schedules (i.e., reliable and less severe patients) are given few slots, some not at their preferred time.  The patients in between (i.e., reliable, but more severe patients) are given many, but undesirable slots.  Minor variations to this logic occur due to the total number of slots available to be assigned.

As shown in Figure \ref{f:I_MSR_VS}, even in settings involving patients with a common time-of-day preference, the population appointment schedule produced by the optimization-based approach is more complex than the one produced by a simple rotation policy, weighing the relative impact of each dimension (time-of-day preference, severity, and reliability) when assigning slots.  The rotation policy simply focuses on diversifying slots and providing equal access to all patients.  For settings involving patients with a common time-of-day preference, such a simple strategy works well, even if patients vary in terms of severity and/or reliability.

\begin{observation}
When patient time-of-day preferences vary, accounting for these differences in the optimization-based approach can lead to improvements of up to 15\% over the simple rotation policy that ignores these differences and treats all patients the same.
\end{observation}

\begin{figure}[htbp]
	\centering
	\includegraphics[width=0.9\textwidth]{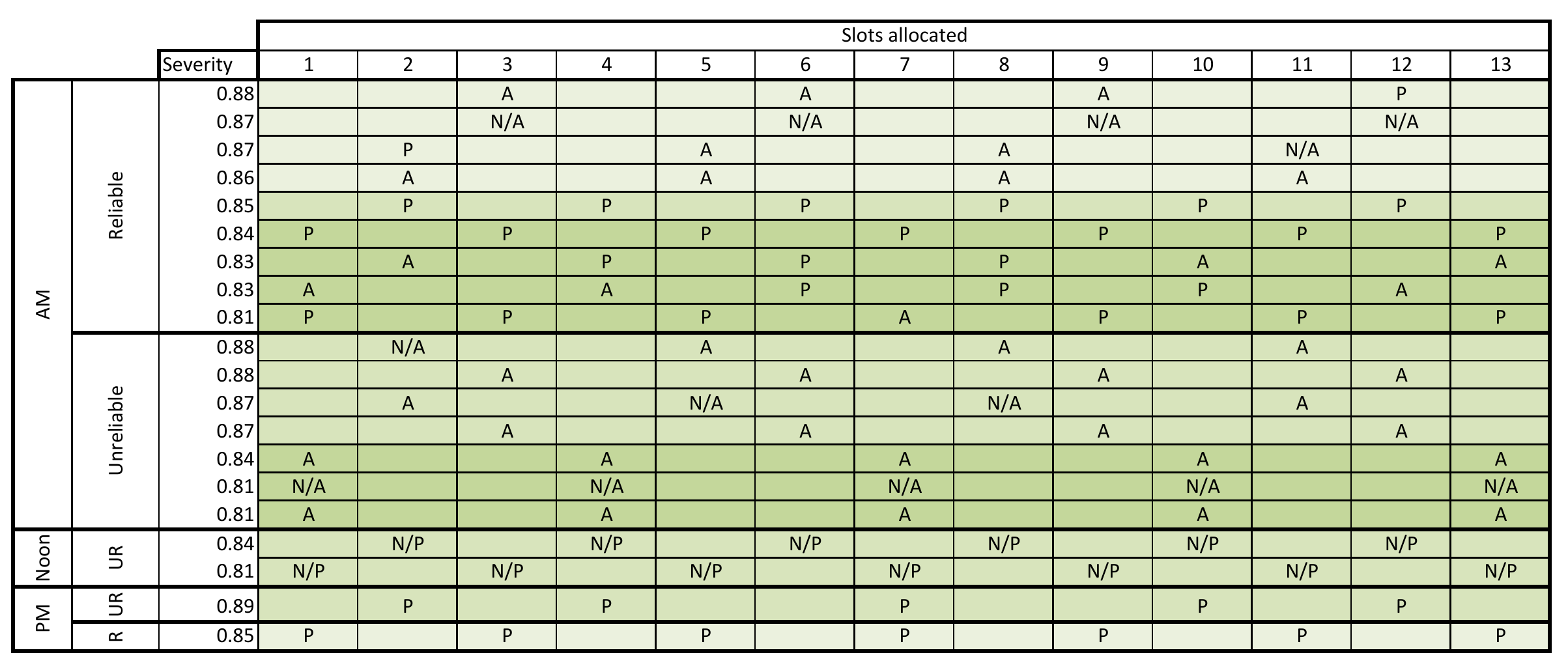}
	\caption{\small Time slots assigned in the optimization-based solution: Profile II: Predominately AM preference; mixed reliability; varying severity: (A) morning (1-3); (P) afternoon (6-8); (N/A) noon/morning (4); (N/P) noon/afternoon.	\label{f:II_MSR_VS}}
\end{figure}

As time-of-day preferences start to vary, the difference in quality of the schedule produced by the simple rotation policy and the schedule produced by the optimization-based approach increases. Even when only 20\% of the patients have a differing time-of-day preference (Profile II), performance gaps between one and two percent become performance gaps between 6.5 and 8.5 percent when time-of-day preferences are strong. In these settings, the schedule produced by the optimization-based approach ensures that, when possible, patients are given their preferred time slots, and, when not possible, a similar logic to what we have seen for the common time-of-day preference setting is employed.
Figure \ref{f:II_MSR_VS} shows the population appointment schedule produced by the optimization-based approach when most patients have a strong preference for morning slots, patients have mixed reliability, and severity levels vary across patients. The noon time slots are further differentiated to account for the fact that some are also desirable for patients with a morning preference and some are  also desirable for patients with an afternoon preference, i.e., a preferred joint noon-AM slot (N/A) or a preferred joint noon-PM slot (N/P).

An examination of the population appointment schedule shows that patients in the two smaller cohorts (with PM and noon time preferences) are assigned their preferred slots (all 13 noon/PM slots are assigned to the patients with a noon preference and 12 of the PM slots are assigned to the patients with an afternoon preference). The allocation of the remaining slots across the patients with a morning preference employs the logic that we have seen before to assign slots to patients with a common time-of-day preference. Only preferred slots are assigned to unreliable patients and more of them if their disease is more severe; if possible, few, but preferred, slots are assigned to reliable patients, but, if not possible, more, but a mix of desirable and undesirable, slots are assigned to reliable patients.

\begin{figure}[htbp]
	\centering
	\includegraphics[width=0.9\textwidth]{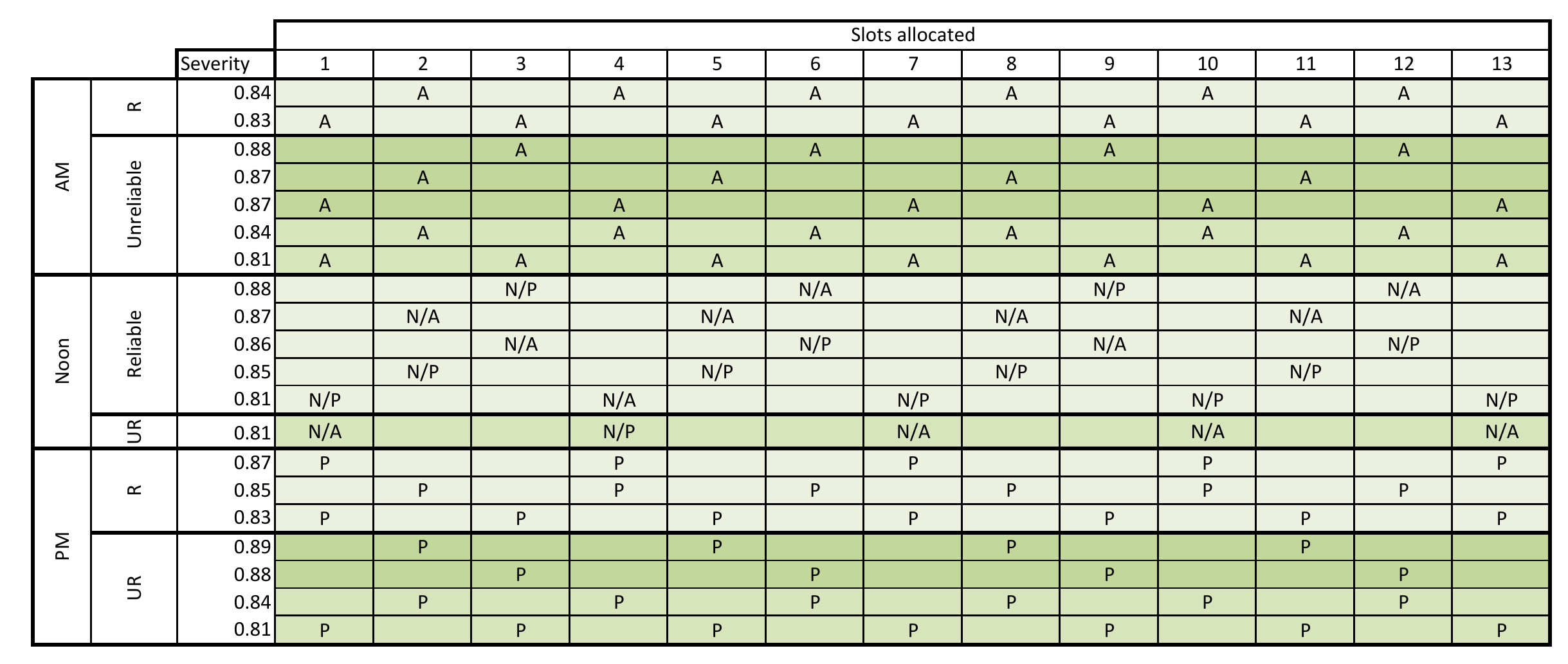}
	\caption{\small Time slots assigned in the optimization-based solution: Profile VI:  Mixed time preference (7 6 7); mixed reliability; varying severity: (A) morning; (P) afternoon; (N) noon. \label{f:VI_MSR_VS}}
\end{figure}

The difference in quality between the rotation policy and the optimization-based approach is largest when the patient population has balanced AM, Noon, PM preferences (i.e., time-of-day patient population profile VI), with performance gaps between 11.5 and 15.5 percent.  In these settings, the schedule produced by the optimization-based approach allocates preferred slots to each group of patients. Of course, such an allocation is naturally imbalanced, because there are fewer noon slots. (The patients with a morning or afternoon preference have three slots per period whereas the patients with a noon preference only have two slots per period.) Consider, Figure \ref{f:VI_MSR_VS}, in which we show the population appointment schedule produced by the optimization-based approach when patients have balanced, but strong, time-of-day preferences, patients have mixed reliability, and severity levels vary across patients.
We see that the average number of visits over the planning horizon is 5.6 for patients with a morning or afternoon preference and only 4.3 for patients with a noon preference.  However, in this specific instance, there is a larger fraction of reliable patients with a noon time preference, compared to the fraction of reliable patients among those with a morning or afternoon preference, and, thus, fewer slots are required to produce effective appointment schedules for the patients  with a noon time preference.  The optimization-based approach ``recognizes'' such instance-specific characteristics and exploits them, whereas the simple rotation policy does not and simply assigns either five or six appointments to patients (with an average of 5.2).

These results suggest that for heterogeneous patient populations, stratifying and scheduling by cohorts can be beneficial.  In the following, we evaluate the ability of simple cohort scheduling policies to capture patient differences and produce high-quality quality schedules.

\vspace{12pt} \noindent \textbf{Analysis of the 1-level cohort scheduling policies}

\begin{observation}
Stratifying patients by time-of-day preference yields significant improvement over the simple rotation policy when the patients in the population have clearly distinguishable time-of-day preferences. Stratifying along other distinguishing factors can lead to low
quality schedules.
\end{observation}

Using a time-based 1-level cohort scheduling policy, which partitions the set of patients into cohorts based on their time-of-day preferences and allocates an appropriate number of slots to each cohort, can significantly increase the quality of the population appointment schedule (compared to the quality of the population appointment schedule produced by the simple rotation policy) as it can avoid assigning undesirable time slots to patients. In Table \ref{t:cohort_perf}, we see that the maximum performance gap for the time-based 1-level cohort strategy is 2.4\% when patients have strong time-of-day preferences and 2.3\% when patients have weak time-of-day preferences. The average performance gap is only 1.1\% when patients have strong preferences and 0.9\% when patients have weak preferences.

Table \ref{t:TP_perf} shows the performance of the time-based 1-level cohort scheduling policy for the different time-of-day preference profiles. Columns 2 - 4 show the fraction of the total number of slots allocated to each cohort, columns 5 - 7 show the average number of slots assigned to a patient for each cohort, columns 8-10 show the fraction of preferred time slots assigned to each cohort, and column 9 shows the average performance gap over all instances with a given profile.
\begin{table}[htbp]
	\centering
	\caption{\small Performance of the time-based 1-level cohort scheduling policies for different time-of-day preference profiles. \label{t:TP_perf}}
	\vspace{6pt}
	\small
	\begin{tabular}{r||rrr||rrr||rrr||r}
		\toprule
		&\multicolumn{3}{c||}{Slot allocation} & \multicolumn{3}{c||}{Slots per patient}  & \multicolumn{3}{c||}{Preferred slots} & \multicolumn{1}{c}{Avg. gap} \\
		\midrule
		Profile & AM    & Noon  & PM    & AM    & Noon  & PM    & AM    & Noon  & PM  & \\ \hline
		II (16-2-2) & $\frac{3}{4}$ & $\frac{1}{8}$ & $\frac{1}{8}$ & 4.9   & 6.5   & 6.5   & $\frac{2}{3}$ & 1   & 1   & 0.9\% \\
		III (10-5-5) & $\frac{1}{2}$ & $\frac{1}{4}$ & $\frac{1}{4}$ & 5.2   & 5.2   & 5.2   & $\frac{2}{3}$ & 1   & 1   & 1.5\% \\
		IV (10-0-10) & $\frac{1}{2}$ & 0 & $\frac{1}{2}$ & 5.2   & -     & 5.2   & 1 &  - & 1   & 0.8\% \\
		V (5-10-5) & $\frac{1}{4}$ & $\frac{1}{2}$ & $\frac{1}{4}$ & 5.2   & 5.2   & 5.2   & 1   & $\frac{1}{2}$ & 1   & 0.7\% \\
		VI (7-6-7) & $\frac{3}{8}$ & $\frac{2}{8}$ & $\frac{3}{8}$ & 5.6   & 4.3   & 5.6   & 1   & 1    & 1   & 1.2\% \\
		\bottomrule
	\end{tabular}%
	\normalsize
\end{table}%

We see that by allocating the chosen fraction of the total number of slots to each of the cohorts, the average number of slots per patient as well as the quality of slots assigned to patients (i.e., whether a preferred or non-preferred slot is assigned) in a cohort is close to ideal. (Note that because there are 13 periods, 8 slots per period, and 20 patients, the average number of slots per patient is 5.2.) As a consequence, the average performance gaps for the different time-of-day preference profiles are small. The population appointment schedules produced by the optimization-based approach are slightly better because they better accommodate differences in reliability and disease severity, and have greater flexibility in the number and spread of visits for a patient over the planning horizon based on the desirability of the slots assigned.

The 1-level cohort scheduling policies that partition patients based on either their disease severity or their reliability (and ignore their time-of-day preferences) result in low-quality population appointment schedules. The average gap for the severity-based 1-level cohort scheduling policy is 13.1\% for strong time preferences and 4.8\% for weak time preferences. The average gap for the reliability-based 1-level cohort scheduling policy is 12.5\% for strong time preferences and 4.9\% for weak time preferences.  We can see that ignoring time-of-day preferences, even when those preferences are weak, leads to poor schedules. When patient time-of-day preferences vary, the performance gaps for the population appointment schedules obtained with the severity-based and reliability-based 1-level cohort strategies are often worse than those obtained with the simple rotation policy.

These insights, of course, have been influenced by the characteristics of the instances in our test set.  However, we are confident that the success of a 1-level cohort strategy depends on two key conditions: (1) the characteristic used to define the cohorts must have a significant impact on the quality of a population appointment schedule, and (2) the cohorts as well as an appropriate allocation of slots to cohorts must be easy to identify.
In the settings considered in our computational study, these two conditions were satisfied for the time-of-day preferences, but not for the disease severity and reliability (in both cases the first condition was not met).

In summary, our computational experiments have demonstrated that accounting for time-of-day preferences of patients can significantly improve population appointment schedule quality. Furthermore, in many situations, relatively simple, but effective time-based cohort scheduling policies can yield populations appointments schedules of similar quality. Creating cohorts based on other patient characteristics did not perform well for the settings considered.

\section{Final Remarks}
\label{s:final}

We have investigated optimization- and cohort-based methods for scheduling appointments for patients in a community-based chronic disease management program with the goal of minimizing the aggregate probability of patients being in an uncontrolled health state. The
optimization-based method explicitly accounts for disease progression since the time of the last appointment and the possibility that patients fail to show up at appointments. Our computational study (1) highlights the considerable impact that time-of-day preferences can have on population health outcomes, (2) demonstrates that simple strategies, i.e., cohort-based scheduling policies, can be effective in reducing no-show rates when the patient population can easily be divided into cohorts with similar time-of-day preferences, and (3) optimization-based methods are preferred and provide better health outcomes when accurate and detailed individual patient information is available. The latter suggests that developing and putting in place processes to gather that data, e.g., by specifically focusing on such data during intake consultations or by including and monitoring operational data within electronic medical records, should be considered as the benefits to population health outcomes of using that information can be substantial.

Our computational study also reveals that the highest quality population appointment schedules carefully tradeoff the visit frequency and the desirability of visit times to control no-show rates. In most situations, intuitive rules-of-thumb, i.e., higher visit frequencies for patients with more severe disease levels, spreading patient visits equally throughout the planning period, and assigning more desirable visit times to patients, perform well and when applied in a straightforward way (as in the cohort scheduling
policies) can substantially improve population health outcomes.

\bibliographystyle{natbib}
\bibliography{pat_assign}


\section*{Acknowledgment}

This research has been supported by grant CMII-0654398 from the
National Science Foundation.  The authors thank Katrese Minor, Dr.
Paul Detjen, and the entire staff at the Mobile C.A.R.E. Foundation,
and Tricia Morphew of the Asthma and Allergy Foundation of America
for their invaluable input on childhood asthma and mobile health
care programs.

\newpage

\section*{Appendix A. Cohort scheduling algorithms}

Given a number of cohorts $C$ and for each cohort $c$, for $c =
1,...,C$, the set of patients $P^c = \{p^c_1, p^c_2, ...,
p^c_{n^c}\}$ and the set of time slots $T^c = \{t^c_1, t^c_2, ...,
t^c_{k^c}\}$, where $\cup^C_{c=1} P^c = \{1,...,P\}$, $\cup^C_{c=1}
T^c = \{1,...,T\}$, and $P^i \cap P^j = \emptyset$ and $T^i \cap T^j
= \emptyset$ for all $i,j = 1,...,C, i \neq j$, Algorithm
\ref{alg:cohort} creates the population appointment schedule.

\vspace{6pt}
\begin{algorithm}[H]
\For{$c \leftarrow 1$ \KwTo $C$}{
    $i \leftarrow 1$ \;
    \For{$k\leftarrow 1$ \KwTo $K$}{
        \For{$j \leftarrow 1$ \KwTo $k^c$}{
            Assign patient $p^c_i$ to time slot $t^c_j$ in period $k$ \;
            \lIf{$i = n^c$}{$i \leftarrow 1$} \lElse{$i \leftarrow i+1$} \;
        }
    }
} \caption{Creating an appointment schedule with a cohort policy.}
\label{alg:cohort}
\end{algorithm}
\vspace{6pt}

The rotation policy naturally introduces diversification in the time
slots assigned to a patient \textit{unless} the number of patients
is a multiple of the number of time slots, because in that case, a
patient will be assigned the same time slot in each of his visits.
When the number of patients is a multiple of the number of time
slots, diversification is accomplished by introducing a
slot-reversing rotation policy as shown in Algorithm \ref{alg:sr}.

\begin{algorithm}[H]
\For{$c \leftarrow 1$ \KwTo $C$}{
    $i \leftarrow 1$ \;
    $direction \leftarrow up $\;
    \For{$k\leftarrow 1$ \KwTo $K$}{
        \For{$j \leftarrow 1$ \KwTo $k^c$}{
            \eIf{direction = up}{
                Assign patient $p^c_i$ to time slot $t^c_j$ in period $k$ \;
            }{
                Assign patient $p^c_i$ to time slot $t^c_{k^c-j}$ in period $k$ \;
            }
            \eIf{$i = n^c$}{
                $i \leftarrow 1$ \;
                \lIf{direction = up}{$direction \leftarrow down$} \lElse{$direction \leftarrow up$} \;
            }
            {
                $i \leftarrow i+1$ \;
            }
        }
    }
} \caption{Creating an appointment schedule with a cohort policy
with slot reversing.} \label{alg:sr}
\end{algorithm}

\clearpage 

\section*{Appendix B. Slot allocations for cohort scheduling algorithms}

\begin{table}[htp]
	\vspace{-0.5cm}
	\caption{Slot allocation for time-of-day cohorts \label{coh:tod}} \vspace{4pt}
	\centering
	\begin{tabular}{r r r| r r r}
		\hline 	
		\multicolumn{3}{c|}{Profile} & \multicolumn{3}{c}{Cohort slot allocation (\% total)} \\
		\cline{4-6}
		AM  &Noon  &PM &Cohort 1 (AM) & Cohort 2 (Noon) & Cohort 3 (PM) \\
		\hline
		20 & 0 & 0 & 1 & 0 & 0 \\
		18 & 2 & 2 & $\frac{3}{4}$ & $\frac{1}{8}$ & $\frac{1}{8}$ \\
		10 & 5 & 5 & $\frac{1}{2}$ & $\frac{1}{4}$ & $\frac{1}{4}$ \\
		10 & 0 & 10 & $\frac{1}{2}$ & 0 & $\frac{1}{2}$ \\
		5 & 10 & 5 & $\frac{1}{4}$ & $\frac{1}{2}$ & $\frac{1}{4}$ \\
		7 & 6 & 7 & $\frac{3}{8}$ & $\frac{2}{8}$ & $\frac{3}{8}$ \\
		\hline 	
	\end{tabular}
\end{table}

\begin{table}[htp]
	\vspace{-0.5cm}
	\caption{Slot allocation for severity cohorts \label{coh:sev}} \vspace{4pt}
	\centering
	\begin{tabular}{r r | r r}
		\hline 	
		\multicolumn{2}{c|}{Profile} & \multicolumn{2}{c}{Cohort slot allocation (\% total)} \\
		\cline{3-4}
		Severe &Mild &Cohort 1 (Severe)& Cohort 2 (Mild) \\
		\hline
		20 & 0 & 1 & 0 \\
		10 & 10 & $\frac{5}{8}$ & $\frac{3}{8}$ \\
		0 & 20 & 0 & 1 \\
		10 ($<0.85$) &10 ($>0.85$) &$\frac{5}{8}$ & $\frac{3}{8}$ \\
		\hline 	
	\end{tabular}
\end{table}

\begin{table}[htp]
\vspace{-0.5cm}
	\caption{Slot allocation for reliability cohorts \label{coh:rel}} \vspace{4pt}
	\centering
	\begin{tabular}{r r | r r}
		\hline 	
		\multicolumn{2}{c|}{Profile} & \multicolumn{2}{c}{Cohort slot allocation (\% total)} \\
		\cline{3-4}
		Reliable &Unreliable & Cohort 1 (Reliable) & Cohort 2 (Unreliable) \\
		\hline
		20 & 0 & 1 & 0 \\
		10 & 10 & $\frac{3}{8}$ & $\frac{5}{8}$ \\
		0 & 20 & 0 & 1 \\
		12 ($<0.9$) &8 ($>0.9$) & $\frac{1}{2}$ & $\frac{1}{2}$ \\
		\hline 	
	\end{tabular}
\end{table}

\clearpage 

\section*{Appendix C. Slot assignments for additional profiles}

\begin{figure}[htp]
    \centering
        \includegraphics[width=0.9\textwidth]{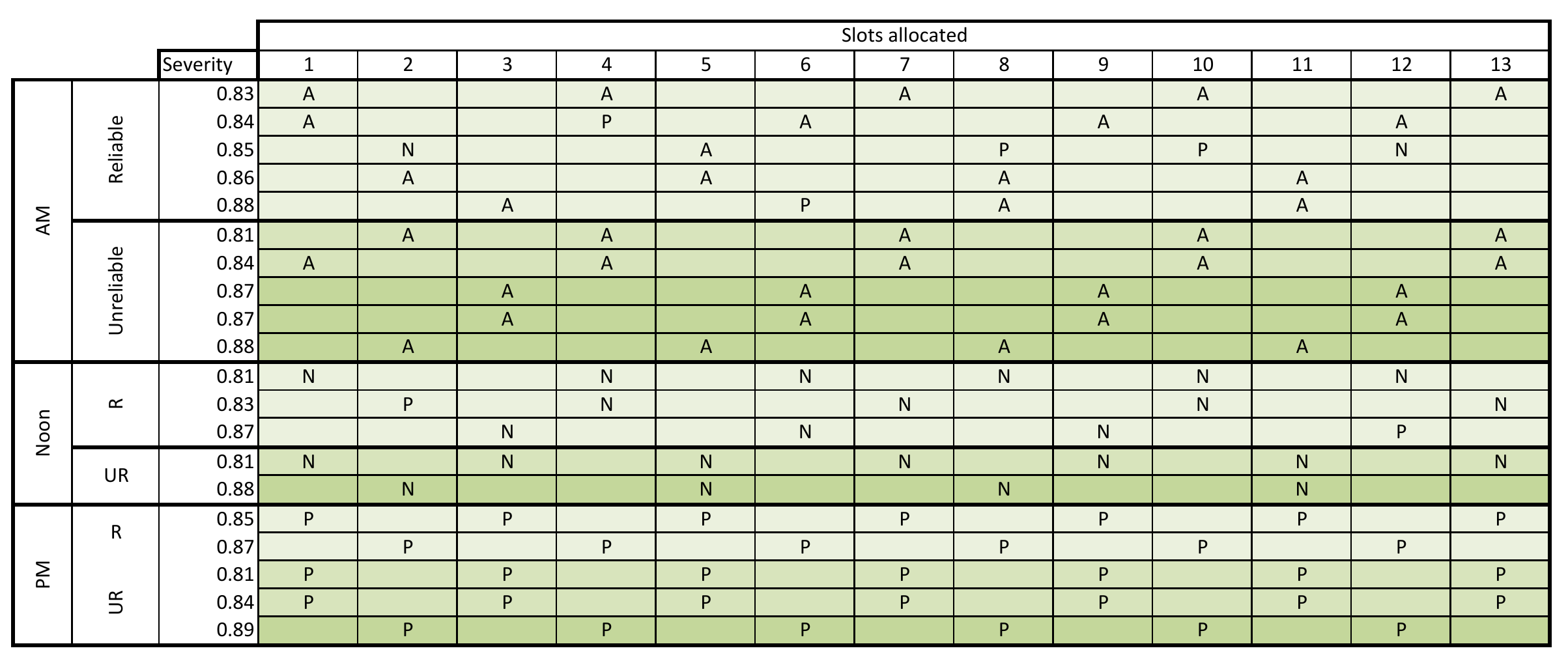}
    \caption{\small Time slots assigned in the optimization-based solution: Profile III: Mixed time preference (10 5 5); mixed reliability; varying severity: (A) morning; (P) afternoon; (N) noon. \label{f:III_MSR_VS}}
\end{figure}

\begin{figure}[htp]
    \centering
        \includegraphics[width=0.9\textwidth]{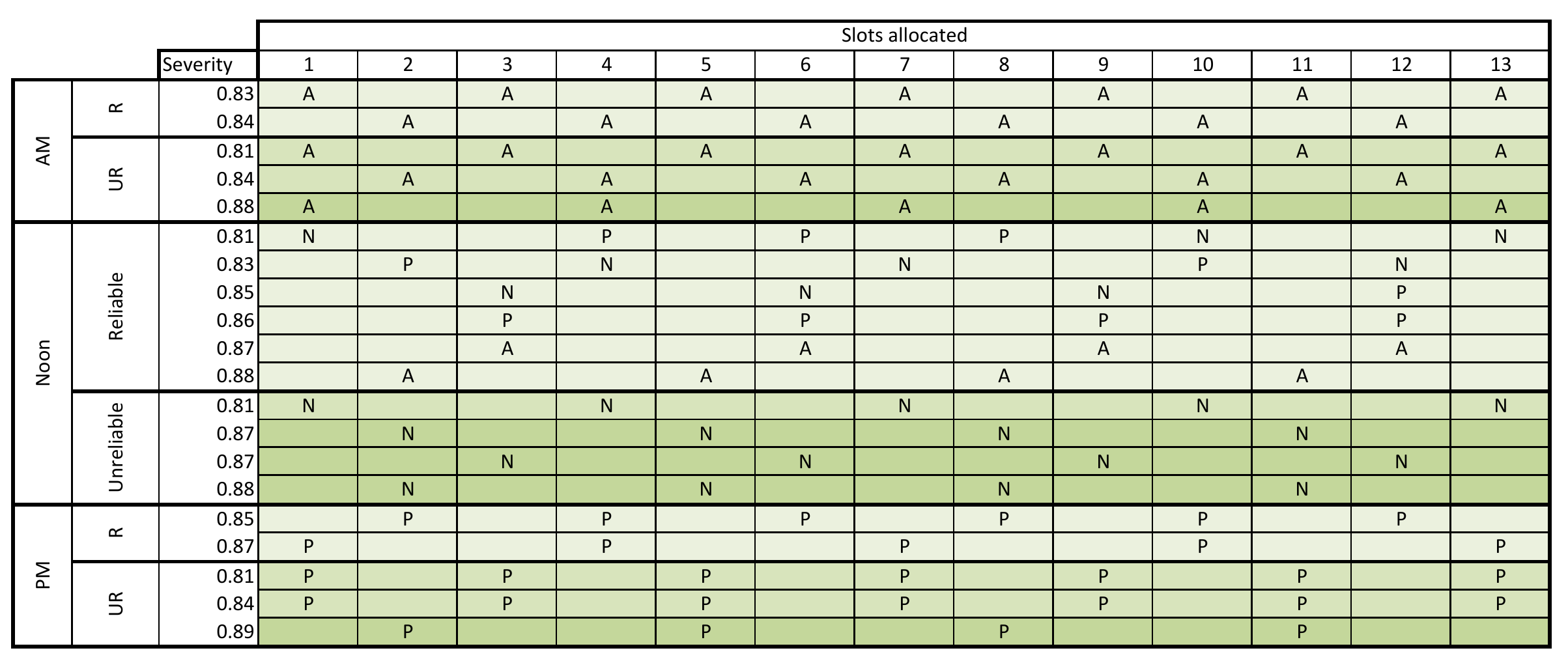}
    \caption{\small Time slots assigned in the optimization-based solution: Profile V: Mixed time preference (5 10 5); mixed reliability; varying severity: (A) morning; (P) afternoon; (N) noon. \label{f:V_MSR_VS}}
\end{figure}

\end{document}